\documentstyle[epsf,12pt]{article}
\catcode`@=11
\renewcommand{\@begintheorem}[2]{\it \trivlist            
      \item[\hskip \labelsep{\bf #1\ #2{\rm :}}]}         
\renewcommand{\@opargbegintheorem}[3]{\it \trivlist       
      \item[\hskip \labelsep{\bf #1\ #2\ {\rm (#3)\/:}}]}
\def\@sect#1#2#3#4#5#6[#7]#8{\ifnum #2>\c@secnumdepth
     \def\@svsec{}\else 
     \refstepcounter{#1}\edef\@svsec{\csname the#1\endcsname{.}\hskip 1em }\fi
     \@tempskipa #5\relax
      \ifdim \@tempskipa>\z@ 
        \begingroup #6\relax
          \@hangfrom{\hskip #3\relax\@svsec}{\interlinepenalty \@M #8\par}
        \endgroup
       \csname #1mark\endcsname{#7}\addcontentsline
         {toc}{#1}{\ifnum #2>\c@secnumdepth \else
                      \protect\numberline{\csname the#1\endcsname}\fi
                    #7}\else
        \def\@svsechd{#6\hskip #3\@svsec #8\csname #1mark\endcsname
                      {#7}\addcontentsline
                           {toc}{#1}{\ifnum #2>\c@secnumdepth \else
                             \protect\numberline{\csname the#1\endcsname}\fi
                       #7}}\fi
     \@xsect{#5}}

\@addtoreset{equation}{section}

\newcommand{\qed}{\hspace*{\fill}\rule{0.25cm}{0.25cm}}

\newtheorem{lemma}{Lemma}[section]
\newtheorem{theorem}[lemma]{Theorem}

\newcounter{rnc}

\newcommand{\rnp}[1]{\setcounter{rnc}{#1}{\rm (\roman{rnc})}}

\newcommand{\Lomit}[1]{} 
\newcommand{\Xomit}[1]{}
\def\c{\widetilde{\rm c}}
\def\p{{\rm P}}
\def\B{{\rm B}}
\def\D{{\cal D}}
\def\P{{\cal P}}
\def\O{{\rm O}}
\def\Z{{\bf Z}}

\def\R{{\bf R}}
\def\0{{\bf 0}}
\def\de{\delta}
\def\phi{\varphi}
\def\back{\backslash}
\def\Push{{\sf Push}}
\def\Scan{{\sf Scan}}
\def\Fix{{\sf Fix}}
\def\SFM{{\sf SFM}}

\begin{document}
\title{A Combinatorial, Strongly Polynomial-Time Algorithm for 
Minimizing Submodular Functions}

\author{Satoru {\sc Iwata}
\thanks{Division of Systems Science, Graduate School of Engineering
Science, Osaka University, Toyonaka, Osaka 560-8531, Japan. 
E-mail: {\tt iwata@sys.es.osaka-u.ac.jp}. A part of this work is 
done while on leave at the Fields Institute, Toronto, Canada. 
Partly supported by Grants-in-Aid for Scientific Research from 
Ministry of Education, Science, Sports, and Culture of Japan.} 
\and
Lisa {\sc Fleischer}
\thanks{Department of Industrial Engineering
and Operations Research, Columbia University, New York, NY 10027, USA.
E-mail: {\tt lisa@ieor.columbia.edu}.  This work done while
on leave at Center for Operations Research and Econometrics, Universit\'e
Catholique de Louvain, Belgium, and at the Fields Institute, Toronto, 
Canada. Partially supported by NSF grants INT-9902663 and EIA-9973858.}
\and
Satoru {\sc Fujishige}
\thanks{Division of Systems Science, 
Graduate School of Engineering Science,  
Osaka University, Toyonaka, Osaka 560-8531, Japan. 
E-mail: {\tt fujishig@sys.es.osaka-u.ac.jp}. 
Partly supported by Grants-in-Aid for Scientific Research from 
Ministry of Education, Science, Sports, and Culture of Japan.}}

\date{July 1999; revised October 1999}

\maketitle

\begin{abstract}
This paper presents the first combinatorial polynomial-time algorithm for 
minimizing submodular set functions, answering an open question posed
in 1981 by Gr\"otschel, Lov\'asz, and Schrijver. 
The algorithm employs a scaling 
scheme that uses a flow in the complete directed graph on the underlying 
set with each arc capacity equal to the scaled parameter. The resulting 
algorithm runs in time bounded by a polynomial in the size of the 
underlying set and the largest length of the function value. The paper 
also presents a strongly polynomial-time version that runs in time bounded 
by a polynomial in the size of the underlying set independent of the 
function value. 
\end{abstract}


\clearpage
\section{Introduction} 
Gr\"otschel, Lov\'asz, and Schrijver~\cite{GLS81} revealed the 
polynomial-time equivalence between the optimization and separation 
problems in combinatorial optimization 
via the ellipsoid method. Since then, many combinatorial 
problems have been shown to be polynomial-time solvable by means of 
their framework. The problem of minimizing submodular (set) functions 
is among these problems. 
Since the ellipsoid method is far from being efficient in practice
and is not combinatorial, efficient combinatorial algorithms for 
submodular function minimization have been desired for a long time. 

A function $f$ on all the subsets of a finite set $V$ is called 
{\em submodular} if it satisfies 
$$f(X)+f(Y)\ge f(X\cup Y)+f(X\cap Y)\quad\quad\forall X, Y\subseteq V.$$
We suppose that $f(\emptyset)=0$ without loss of generality throughout 
this paper. 

Submodular functions arise in various branches of mathematical engineering
such as combinatorial optimization and information theory.  There
are also close connections between submodularity and 
convexity~\cite{Fuji91,Lovasz83}.
\Xomit{
See Lov\'asz~\cite{Lovasz83} and Fujishige~\cite{Fuji91} for 
fundamental results about submodular functions and for close connections 
to convexity.  
}
Examples include the matroid rank function, the cut capacity function, 
and the entropy function. 
In each of these and other applications, 
the {\em base polyhedron} defined by 
\begin{equation}
\label{eq:base}
\B(f)=\{x\mid x\in \R^V,\, x(V)=f(V),\,
\forall X\subseteq V:\ x(X)\le f(X)\}
\end{equation}
often plays an important role, 
where $x(X)=\sum_{v\in X}x(v)$ for any $X\subseteq V$.

Linear optimization problems over base polyhedra are efficiently 
solvable by the greedy algorithm of Edmonds~\cite{Edm70}. 
Thus Gr\"otschel, Lov\'asz, and Schrijver \cite{GLS81} assert that 
the submodular function minimization, which is equivalent to the 
separation problem, is solvable in polynomial time by the ellipsoid 
method. Later, they also devise a strongly polynomial-time algorithm 
within their framework using the ellipsoid method \cite{GLS88}.   

A first step towards a combinatorial strongly polynomial-time 
algorithm was taken by Cunningham \cite{Cunningham84,Cunningham85}, 
who devised a strongly polynomial-time algorithm for testing membership 
in matroid polyhedra as well as a pseudopolynomial-time algorithm for 
minimizing submodular functions. Recently, Narayanan \cite{Narayanan95} 
improved the running time bounds of these combinatorial algorithms by 
a rounding technique. Based on the minimum-norm base characterization 
of minimizers due to Fujishige \cite{Fuji80,Fuji84}, 
Sohoni~\cite{Sohoni92} gave another combinatorial pseudopolynomial-time 
algorithm for submodular function minimization. 

For the problem of minimizing a symmetric submodular function over 
proper nonempty subsets, Queyranne \cite{Queyranne98} presented 
a combinatorial strongly polynomial-time algorithm, extending 
the undirected minimum cut algorithm of Nagamochi and Ibaraki 
\cite{Nagamochi+Ibaraki92}. 

In this paper, we present a combinatorial polynomial-time algorithm 
for submodular function minimization.  Our algorithm uses an augmenting
path approach with reference to a convex combination of extreme points of
the base polyhedron.
Such an approach was first introduced by Cunningham for minimizing 
submodular functions that arise from the separation problem 
for matroid polyhedra~\cite{Cunningham84}.
This was adapted for general submodular function minimization
by Bixby, Cunningham, and Topkis~\cite{Bixby+Cunningham+Topkis} 
and improved by Cunningham~\cite{Cunningham85} to obtain 
a pseudopolynomial-time algorithm.  
\Xomit{
We succeed Cunningham's approach 
of finding an augmenting path with reference to a convex combination of 
extreme bases. 
}


A fundamental tool in these algorithms is to move from one extreme
point of the base polyhedron to an adjacent extreme point 
via an exchange operation that
increases one coordinate and decreases another coordinate by the
same quantity. 
This quantity is called the exchange capacity.
These previous methods maintain a directed graph on the underlying
set that represents the possible exchange operations.
They are inefficient since the lower bound on the size of each augmentation 
is too small. In traditional network flow problems, it is possible to 
surmount this difficulty by augmenting only on paths of sufficiently large 
capacity~\cite{Edmonds+Karp72}.  However, it has proved difficult to 
adapt this scaling approach to work in the setting of 
submodular function minimization, mainly because the amount of
augmentation is determined by exchange capacities multiplied by 
the convex combination coefficients.  These coefficients can 
be as small as the 
reciprocal of the maximum absolute value of the submodular function. 

\Xomit{in part because the set of arcs used in these algorithms is
the set of arcs in the Hasse diagrams of the posets corresponding
to tight sets of extreme bases,
and is not the full
set of arcs with positive exchange capacity, and in part because
the capacity of a path is not determined solely by the separate
capacities of each arc in the path.}

To overcome this difficulty, we augment the directed graph corresponding
to allowable exchanges
with the complete directed graph on the underlying set, 
letting the capacity of this additional arc set depend directly 
on our scaling parameter. This technique was first introduced 
by Iwata~\cite{Iwata97}, who used it to develop the first 
polynomial-time capacity-scaling algorithm for the submodular 
flow problem of Edmonds and Giles~\cite{Edmonds+Giles77}.
This algorithm was later refined by 
Fleischer, Iwata, and McCormick~\cite{FIM99} into one of the
fastest algorithms for submodular flow.  Our work in this paper 
builds on ideas in this latter paper to develop a capacity-scaling,
augmenting-path algorithm for submodular function minimization.
The running time 
of the resulting algorithm is weakly polynomial, i.e., bounded by 
a polynomial in the size of the underlying set and the largest length 
of the function value. 
Even under the similarity assumption that 
the largest function value is bounded by a polynomial in
the size of the underlying set, our algorithm is faster than
the best previous combinatorial, pseudopolynomial-time
algorithm~\cite{Cunningham85}.

\Xomit{
We adopt a scaling scheme that utilizes flows in the 
complete graph on the underlying set with each arc capacity equal to 
the scaled parameter. Use of such graphs in scaling submodular 
functions was originally made by Iwata \cite{Iwata97} to extend the 
Edmonds--Karp minimum-cost flow algorithm \cite{Edmonds+Karp72} to 
the submodular flow problem of Edmonds and Giles \cite{Edmonds+Giles77}. 
Fleischer, Iwata, and McCormick \cite{FIM99} have substantially improved 
the running time of the capacity scaling algorithm by handling the 
additional flow more effectively. The present paper applies this refined 
scaling technique to submodular function minimization.}

We then modify our scaling algorithm to run in strongly polynomial time, 
i.e., in time bounded by a polynomial in the size of the 
underlying set, independently of the largest length of the function value. 
To make a weakly 
polynomial-time algorithm run in strongly polynomial time, 
Frank and Tardos \cite{Frank+Tardos87} developed a generic preprocessing 
technique that is applicable to a fairly wide class of combinatorial 
optimization problems including the submodular flow problem and testing
membership in matroid polyhedra. However, this framework does not apply to 
submodular function minimization. Instead, we devise a combinatorial 
algorithm that repeatedly detects an element that belongs to every 
minimizer or an ordered pair of elements with the property that if the first 
belongs to a minimizer then the second does.

There are some practical problems, in dynamic flows
\cite{Hoppe+Tardos95}, facility location \cite{Tamir93}, and 
multi-terminal source coding \cite{Fuji78,Han79}, 
where the polynomial-time solvability relies on a submodular 
function minimization routine. 
Goemans and Ramakrishnan \cite{Goemans+Ramakrishnan95} discussed 
a class of submodular function minimization problems 
over restricted families of subsets. Their solution is combinatorial 
modulo an oracle for submodular function minimization on distributive 
lattices. Our algorithm can be used to provide combinatorial, strongly 
polynomial-time algorithms for these problems.

This paper is organized as follows. Section~\ref{sec:prl} provides 
preliminaries on submodular functions. Section~\ref{sec:sa} 
presents a scaling algorithm for submodular function minimization, 
which runs in weakly polynomial time.  
Section~\ref{sec:spa} is devoted to the strongly polynomial-time 
algorithm. Finally, we discuss extensions in Section~\ref{sec:cr}.

\section{Preliminaries}
\label{sec:prl}
We denote by $\Z$ and $\R$ the set of integers and the set of 
reals, respectively. Let $V$ be a finite nonempty set of cardinality $|V|=n$.  
For a vector $x\in\R^V$ we define a modular function $x:2^V\to\R$ by 
$x(X)=\sum_{v\in V} x(v)$. 
For each $u\in V$, we denote by $\chi_u$ the unit vector in $\R^V$ 
such that $\chi_u(v)=1$ if $v=u$ and $=0$ otherwise.

\Xomit{
A function 
$f: 2^V \to\R$ is said to be {\it submodular} if it satisfies 
$$f(X)+f(Y)\ge f(X\cup Y)+f(X\cap Y) \quad\quad  (X, Y\subseteq V).$$
We suppose that $f(\emptyset)=0$ without loss of generality throughout 
this paper. 
}
\Xomit{
We define the {\em submodular polyhedron} $\p(f)$ and the 
{\it base polyhedron} $\B(f)$ associated with the submodular 
function $f$ by 
\begin{eqnarray*}
\p(f)&=&\{x\mid x\in \R^V,\, \forall X\subseteq V:\ x(X)\le f(X)\}, \\
\B(f)&=&\{x\mid x\in \p(f),\, x(V)=f(V)\}.
\end{eqnarray*}
}
Given a submodular function $f$ with $f(\emptyset)=0$ and its
associated base polyhedron $\B(f)$ as defined in (\ref{eq:base}),
we call a vector $x\in\B(f)$ a {\it base}. 
An extreme point of $\B(f)$ is called an {\it extreme base}. 
A fundamental step in submodular function minimization algorithms
is to move from one base $x$ to another base $x'$
via an {\em exchange operation} that increases one coordinate while
decreasing another coordinate by the same amount.  The maximum amount of
increase that ensures $x'\in\B(f)$ is called the exchange capacity.  
More precisely, for any base $x\in\B(f)$
and any distinct $u,v\in V$ the {\it exchange capacity} is
\begin{equation}
\label{eq:ecmax}
\c(x,u,v)=
\max\{\alpha\mid\alpha\in{\bf R}, x+\alpha(\chi_u - \chi_v)\in \B(f)\}.
\end{equation}
The exchange capacity can also be expressed as 
\begin{equation}
\label{eq:ecmin}
\c(x,u,v)=\min\{f(X)-x(X)\mid u\in X \subseteq V\back\{v\}\}.
\end{equation}
In general, computing $\c(x,u,v)$ is as hard as submodular function
minimization, even when $x$ is an extreme base.  However,
if $x$ is an extreme base, then for special pairs of vertices $u$ and $v$,
the exchange capacity $\c(x,u,v)$ can be computed with one function 
evaluation as follows.  

%


Let $L=(v_1,v_2,\cdots,v_n)$ be a linear ordering of $V$. For any 
$k\in\{1,2,\cdots,n\}$, we define $L(v_k)=\{v_1,v_2,\cdots,v_k\}$. 
Given such a linear ordering, the greedy algorithm of Edmonds~\cite{Edm70} 
computes 
\begin{equation}\label{eq:greedy0}
  y(v_i)=f(L(v_i))-f(L(v_{i-1}))\quad\quad (i=1,2,\cdots,n),  
\end{equation} 
where $L(v_0)=\emptyset$. The resulting vector $y\in\R^V$ is an extreme 
base $y\in\B(f)$. Conversely, any extreme base can be generated by applying 
the greedy algorithm to an appropriate linear ordering.  
Note that a linear ordering 
$L=(v_1,v_2,\cdots,v_n)$ generates an extreme base $y$
if and only if $y(L(v_i))=f(L(v_i))$ for $i=1,2,\cdots,n$.
For any base $y\in\B(f)$, a set $X\subseteq V$ is called {\it $y$-tight} if
$y(X)=f(X)$.
A pair $(u,v)$ is called {\em eligible} for $y$ if $u$ immediately succeeds 
$v$ in some linear ordering that generates $y$. The following lemma enables 
us to compute an exchange capacity $\c(y,u,v)$ if $(u,v)$ is eligible 
for $y$. 

\begin{lemma}\label{lemma:greedy}
Let $L$ be a linear ordering of $V$ that generates an extreme base 
$y\in\B(f)$. 
Let $L'$ be the linear ordering obtained by interchanging $u$ and $v$ that 
are consecutive in $L$. Then the extreme base $y'$ generated by $L'$ 
satisfies   
\begin{equation}\label{eq:greedy1}
   y'=y+\beta(\chi_u-\chi_v)
\end{equation}
with 
\begin{equation}\label{eq:greedy2}
  \beta=f(L(u)\back\{v\})-f(L(u))+y(v). 
\end{equation} 
Moreover, we have $\c(y,u,v)=\beta$. 
\end{lemma}
\begin{proof}
Equations (\ref{eq:greedy1}) and (\ref{eq:greedy2}) follow from 
the greedy algorithm (see (\ref{eq:greedy0})).  By the definition 
(\ref{eq:ecmax}) of the exchange capacity, we have $\beta\leq\c(y,u,v)$. 
Since $y(L(u))=f(L(u))$, it follows from (\ref{eq:ecmin}) and 
(\ref{eq:greedy2}) that $\beta\geq\c(y,u,v)$. Thus we obtain 
$\beta=\c(y,u,v)$. 
\end{proof}

We will use Lemma~\ref{lemma:greedy} to transform one extreme base into another
and to update the corresponding linear ordering.


For any vector $x\in\R^V$, we denote by $x^-$ the vector in $\R^V$ 
defined by $x^-(v)=\min\{0,x(v)\}$ for $v\in V$. The following 
fundamental lemma easily follows from a theorem of Edmonds~\cite{Edm70} 
on the vector reduction of polymatroids 
(see \cite[Corollaries 3.4 and 3.5]{Fuji91}). 
\begin{lemma}\label{lem:Edm}
For a submodular function $f: 2^V\to\R$ we have 
$$\max\{x^-(V)\mid x\in\B(f)\}=\min\{f(X)\mid X\subseteq V\}.$$
If $f$ is integer-valued, then the maximizer $x$ can be chosen from among 
integral bases. 
\qed
\end{lemma}

We will not use the integrality property indicated in the latter half of
this lemma. Lemma \ref{lem:Edm} shows a min-max relation of strong duality.
A weak duality is described as follows: For any base $x\in{\rm B}(f)$ and
any $X\subseteq V$ we have $x^-(V)\le f(X)$. We call the difference 
$f(X)-x^-(V)$ a {\it duality gap}. Note that, if $f$ is integer-valued 
and the duality gap $f(X)-x^-(V)$ is less than one for some $x\in\B(f)$ and 
$X\subseteq V$, then $X$ minimizes $f$.


\section{A Scaling Algorithm}
\label{sec:sa}
In this section, we describe a combinatorial algorithm for minimizing an 
integer-valued submodular function $f: 2^V\to\Z$ with $f(\emptyset)=0$. 
We assume an evaluation oracle 
for the function value of $f$. Let $M$ denote an upper bound on $|f(X)|$ 
among $X\subseteq V$. 
Note that we can easily compute $M$ by $\O(n)$
calls for the evaluation oracle as follows. Let $y$ be an extreme 
base generated by a linear ordering $L$. For any $X\subseteq V$, we have 
$y^-(V)\leq y(X)\leq f(X)\leq\displaystyle\sum_{v\in V}\max\{0,f(\{v\})\}$. 
Thus we obtain 
$M=\max\{|y^-(V)|,\displaystyle\sum_{v\in V}\max\{0,f(\{v\})\}\}$. 

\subsection{Algorithm Outline}

As indicated earlier, our algorithm uses an augmenting path approach
to submodular function 
minimization~\cite{Bixby+Cunningham+Topkis,Cunningham84,Cunningham85}.
As with previous algorithms, we
maintain a base $x\in \B(f)$ as a convex combination of extreme bases
$y_i\in \B(f)$ indexed by $i\in I$, so that $x=\sum_{i\in I} \lambda_i y_i$.  
Roughly speaking, these previous algorithms use a directed graph
with the arc set defined by the pairs of vertices that are eligible for 
some $y_i$, $i\in I$.
They seek to increase $x^-(V)$ 
by performing exchange operations along
a path of arcs from vertices $s$ with $x(s) < 0$ to vertices $t$ with 
$x(t) >0$.
The algorithms stop with an optimal $x$ when there are no more
augmenting paths. The corresponding minimizer $X$ is 
determined by the set of vertices reachable from vertices $s$ with 
$x(s) < 0$.

To adapt this procedure to a scaling framework, we
use 
a complete directed graph on $V$ with arc capacities 
that depend directly on our scaling parameter $\delta$, 
an idea first introduced for submodular flows in~\cite{Iwata97}. 
\Xomit{
$\phi$ is a flow 
(a function on the arc set of the complete directed
graph obeying capacity constraints that depend on $\delta$) with
{\em boundary} 
\begin{equation}
\partial \phi := \sum_{u\in V}
\end{equation}
}
Let $\phi: V\times V\to\R$ be {\em skew-symmetric}, i.e.,  
$\phi(u,v)+\phi(v,u)=0$ for $u,v\in V$, and {\em $\delta$-feasible} 
in that it satisfies capacity constraints 
$-\delta\leq\phi(u,v)\leq\delta$ for every $u,v\in V$. 
The function $\phi$ can be regarded as a flow in the complete directed 
graph $G=(V,E)$ with the vertex set $V$ and the arc set $E=V\times V$.  
The {\it boundary} 
$\partial\phi: V\to\R$ of $\phi$ is defined by
\begin{equation}
\partial\phi(v)=\sum_{u\in V}\phi(u,v)\quad\quad (v\in V).
\end{equation}
Instead of trying to maximize $x^-(V)$ directly, we define $z=x-\partial\phi$. 
Our algorithm seeks to maximize $z^-(V)$ and thereby increases $x^-(V)$. 

We also maintain linear orderings $L_i$ for $i\in I$ and extreme bases 
$y_i$ generated by them.
\Xomit{
We also do not maintain $A(y)$ explicitly for each extreme base $y$.
Instead, we maintain a linear ordering $\leq_{y}$ that defines a linear extension 
of the poset ${\cal P}(y)=(V,\preceq_y)$, i.e., $v\preceq_y u$ implies 
$v\leq_y u$. If $v\neq u$, $v\leq_y u$, $v\in J_y(u)$, and there is no 
$w\in V$ with $v <_y w <_y u$, then we know that $(u,v)\in A(y)$.  This
is all we will need.
} 
We start with an arbitrary linear ordering $L$ on $V$ and the extreme base
$x\in\B(f)$ generated by $L$.
In addition, we start with the zero flow $\phi={\bf 0}$.  
Thus, initially $z^-(V) = x^-(V) \geq -nM$.
We seek to increase $z^-(V)$, and
in doing so, obtain improvements in $x^-(V)$, via the $\delta$-feasibilty
of $\phi$.

The algorithm consists of scaling phases with a positive parameter 
$\delta$. It starts with $\delta=M$, cuts $\delta$ in half at the 
beginning of each scaling phase, and ends with 
$\delta<1/n^2$. 
\Xomit{
Each $\delta$-scaling phase keeps $\phi$ that satisfies
the capacity constraint 
$-\delta\leq\phi(u,v)\leq\delta$ for every $u,v\in V$. 
We call such a flow {\it $\delta$-feasible}.
}
Each $\delta$-scaling phase maintains a $\delta$-feasible flow $\phi$,
and uses the {\em residual 
graph} $G(\phi)=(V,E(\phi))$ with the arc set 
\begin{equation}
E(\phi)=\{(u,v)\ |\ u,v\in V,\,u\neq v,\,\phi(u,v)\le 0\}.
\end{equation}
Intuitively, $E(\phi)$ consists of the arcs through which we can 
augment the flow $\phi$ by $\delta$ without violating the capacity
constraints. 

\Xomit{
We first give an overview of how the $\delta$-scaling phase works. 
While the goal of the algorithm is to increase $x^-(V)$, we pay 
attention to $z=x-\partial\phi$ rather than $x$ itself. The 
$\delta$-scaling phase primarily aims at increasing $z^-(V)$, 
which brings about an improvement of $x^-(V)$ as well by the 
$\delta$-feasibility of $\phi$. 
}
A $\delta$-scaling phase starts by preprocessing $\phi$ to make it
$\delta$-feasible, and then repeatedly searches to send flow
along augmenting paths in $G(\phi)$ from
$S:= 
\{ v \mid v\in V,\,z(v) \leq -\delta\}$
to $T:= 
\{ v\mid v\in V,\,z(v) \geq \delta\}$.
Such a directed path is called a {\it $\delta$-augmenting path}.

If there are no $\delta$-augmenting paths, then the algorithm checks whether
there is a pair $(u,v)$ of vertices 
such that $u$ is reachable from $S$ by a path with residual capacity 
$\geq \de$, $v$ is not,
and $u$ immediately succeeds $v$ in a linear ordering that generates $y_i$ for
some $i\in I$.  We perform the appropriate exchange operation, and
modify $\phi$ by creating residual capacity on $(u,v)$ 
so that $z=x-\partial\phi$ is invariant.  This
operation may increase the set of vertices reachable from $S$ on
paths of residual capacity $\geq \de$.
\Lomit{
If any such pair exists, we transform $y_i$ by 
$y_i:=y_i+\c(y_i,u,v)(\chi_u-\chi_v)$ and $x$ by 
$x:=x+\alpha(\chi_u-\chi_v)$ with 
$\alpha=\min\{\delta,\lambda_i\c(y_i,u,v)\}$. 
We also change $\phi$ 
If in addition 
$\lambda_i\c(y_i,u,v)>\delta$, we add a new index $k$ to $I$, with $y_k$ 
being equal to the previous $y_i$, to express $x$ as a convex combination. 
In this case, $v$ becomes reachable from $S$ in $G(\phi)$.
}
Once a $\delta$-augmenting path is 
found, the algorithm augments the flow $\phi$ by $\delta$ through 
the path without changing $x$. As a consequence, $z^-(V)$ increases 
by $\delta$ in one iteration. 
This is an extension of a technique for handling exchange
capacity arcs in submodular flows first developed in~\cite{FIM99}.

\subsection{Algorithm Details}

\begin{figure}
\framebox[\textwidth]{\parbox{6in}{
\begin{tabbing}
****\=***\=***\=***\=***\= \hspace{2.5in} \=  \kill
$\SFM(f)$: \\
\> \\
{\tt Input:}\ $f:2^V\rightarrow \Z$ \\
{\tt Output:} $X\subseteq V$ minimizing $f$ \\ 
\> \\
{\bf Initialization:} \\
\> $L_i\gets$ an linear ordering on $V$ \\ 
\> $x\gets$ an extreme base in $\B(f)$ generated by $L_i$ \\ 
\> $I\gets\{i\}$, $y_i\gets x$, $\lambda_i\gets 1$, \\ 
\> $\phi\gets\0$, \\  
\> $\delta\gets M$ \\
{\bf While} $\delta \geq 1/n^2$ {\bf do} \\
\> $\delta\gets\delta/2$ \\
\> {\bf For} $(u,v)\in E$ {\bf do} \\ 
\> \> {\bf If} $\phi(u,v)>\delta$ {\bf then} $\phi(u,v)\gets \delta$ \\
\> \> {\bf If} $\phi(u,v)<-\delta$ {\bf then} $\phi(u,v)\gets -\delta$ \\ 
\> $S\gets\{v\mid x(v)\leq \partial\phi(v)-\delta\}$ \\
\> $T\gets\{v\mid x(v)\geq \partial\phi(v)+\delta\}$ \\ 
\> $W\gets$ the set of vertices reachable from $S$ in $G(\phi)$ \\
\> $Z\gets$ the set of active pairs $(i,v)$ of $i\in I$ and $v\in V$ \\
\Lomit{
\> {\bf While} $S\neq\emptyset$, $T\neq\emptyset$ and 
$\Delta^+W\neq\emptyset$ {\bf do}\\
\> \> {\bf While} $W\cap T=\emptyset$ and $\Delta^+W\neq\emptyset$
{\bf do} \\ 
}
\> {\bf While} $S\neq\emptyset$, $T\neq\emptyset$ and 
$Z \neq \emptyset$ {\bf do}, \\
\> \> {\bf While} $W\cap T=\emptyset$ and 
$Z \neq \emptyset$ {\bf do}, \\
\> \> \> Find an active pair $(i,v)\in Z$. \\
\Lomit{
\> \> \> $Z\gets \{v\mid v\in V\back W,\,(u,v)\in A(y_i)\}$ \\ 
\> \> \> {\bf Repeat} \\
\> \> \> \> Find a vertex $v\in Z$ and $\Push(i,u,v)$. \\
\> \> \> \> Update $W$ and $Z$. \\ 
\> \> \> {\bf until} $Z = \emptyset$ or $|W|$ increases. \\
}
\> \> \> Let $u$ be the vertex succeeding $v$ in $L_i$. \\
\> \> \> Apply $\Push(i,u,v)$. \\  
\> \> \> Update $W$ and $Z$. \\
\> \> {\bf If} $W\cap T\neq\emptyset$ {\bf then} \\
\> \> \> Let $P$ be a directed path from $S$ to $T$ in $G(\phi)$. \\
\> \> \> {\bf For} $(u,v)\in P$ {\bf do} 
$\phi(u,v)\gets\phi(u,v)+\delta$, $\phi(v,u)\gets\phi(v,u)-\delta$ \\
\> \> \> Update $S$, $T$, $W$, and $Z$. \\
\> \> Express $x$ as $x=\sum_{i\in I}\lambda_iy_i$ 
by possibly smaller affinely independent \\
\> \> \> \> subset $I$ and positive coefficients 
$\lambda_i>0$ for $i\in I$. \\
{\bf If} $S=\emptyset$ {\bf then} $X=\emptyset$ 
{\bf else if} $T=\emptyset$ {\bf then} $X=V$ 
{\bf else} $X=W$ \\
{\bf End}.
\end{tabbing}
}}
\caption{A scaling algorithm for submodular function minimization.}
\label{fig:sfm}
\end{figure}

We now describe the scaling algorithm more precisely. Figure \ref{fig:sfm} 
provides a formal description.

At the beginning of the $\delta$-scaling phase, after $\delta$ is cut in 
half, the current flow $\phi$ is $2\delta$-feasible. Then the algorithm 
modifies each $\phi(u,v)$ to the nearest value within the interval 
$[-\delta,\delta]$ to make it $\delta$-feasible.
This may decrease $z^-(V)$ for $z=x-\partial\phi$ by at most 
${n\choose 2}\delta$. The rest of the $\delta$-scaling phase aims at 
increasing $z^-(V)$ by augmenting flow along $\delta$-augmenting paths. 

\Xomit{
Given a base $x=\sum_{i\in I}\lambda_iy_i$ and a $\delta$-feasible flow 
$\phi$, let
$S=\{v\mid x(v)\leq \partial\phi(v)-\delta\}$ and 
$T=\{v\mid x(v)\geq \partial\phi(v)+\delta\}$. 
A directed path from $S$ 
to $T$ in $G(\phi)$ is called a {\it $\delta$-augmenting path}.

}
Let $W$ denote the set of vertices reachable by directed paths from $S$ 
in $G(\phi)$. 
For each $i\in I$ we keep a linear ordering $L_i$ that generates $y_i$. 
We call a vertex $v\in V\back W$ {\it active} if $v$ is the 
last vertex in $L_i$ among vertices in $V\back W$
that satisfies $W\back L_i(v)\neq\emptyset$. 
If $v$ is active in $L_i$, we call $(i,v)$ an {\it active pair}. 
We denote by $Z$ the set of the current active pairs. 

\Lomit{
Then $U_i$ is empty if and only if no arc in $H(y_i)$ leaves $W$. 
}

If $W\cap T=\emptyset$, there is no $\delta$-augmenting path in
$G(\phi)$. Then, as long as there is an active pair $(i,v)$, i.e., 
$Z\neq\emptyset$, 
\Lomit{
 Hasse diagram $H(y_i)$ for some 
$i\in I$ has an arc leaving $W$, 
}
the algorithm repeatedly picks an active pair $(i,v)\in Z$ and 
applies $\Push(i,u,v)$ to $u$ that succeeds $v$ in $L_i$. 
Note that $v$ active implies that $u\in W$.
\Lomit{
applying the operation 
$\Push(i,u,v)$ to each arc $(u,v)\in A(y_i)$ with $v\in V\back W$.
}

The operation $\Push(i,u,v)$ 
starts with reducing the flow through $(u,v)$ by 
$\alpha=\min\{\delta,\lambda_i\c(y_i,u,v)\}$. 
The boundary $\partial\phi$ moves to $\partial\phi+\alpha(\chi_u-\chi_v)$. 
The operation $\Push(i,u,v)$ is called 
{\it saturating} if $\alpha=\lambda_i\c(y_i,u,v)$. Otherwise, it is called 
{\it nonsaturating}. 
A nonsaturating $\Push(i,u,v)$ adds to $I$ a new index $k$ with $y_k:=y_i$,
$\lambda_k:=\lambda_i-\alpha/\c(y_i,u,v)$, and $L_k:=L_i$. Whether the 
$\Push(i,u,v)$ is saturating or not, it updates
$y_i$ as $y_i:=y_i+\c(y_i,u,v)(\chi_u-\chi_v)$, 
$\lambda_i:=\alpha/\c(y_i,u,v)$ if $\c(y_i,u,v)>0$, 
and $L_i$ by interchanging $u$ and $v$. 
\Xomit{
A saturating $\Push(i,u,v)$ updates 
$y_i$ as $y_i:=y_i+\c(y_i,u,v)(\chi_u-\chi_v)$ while nonsaturating one 
adds to $I$ a new index $k$ with $y_k:=y_i+\c(y_i,u,v)(\chi_u-\chi_v)$ and 
$\lambda_k=\alpha/\c(y_i,u,v)$ and updates $\lambda_i$ as 
$\lambda_i:=\lambda_i-\alpha/\c(y_i,u,v)$. 
In either case,
a consistent linear ordering for the new extreme base is the linear ordering
for the original extreme base with $u$ and $v$ swapped;
and 
}
Then the current base $x$ moves to $x+\alpha(\chi_u-\chi_v)$. 
Thus $z=x-\partial\phi$ is invariant. 

\begin{figure}
\begin{center}
\framebox[3in]{\parbox{2in}{
$\Push(i,u,v)$: 
\begin{tabbing}
****\=***\=***\=***\=***\= \hspace{2.5in} \=  \kill
$\alpha\gets\min\{\delta,\,\lambda_i\c(y_i,u,v)\}$ \\
$\phi(u,v)\gets\phi(u,v)-\alpha$ \\
$\phi(v,u)\gets\phi(v,u)+\alpha$ \\
{\bf If} $\alpha<\lambda_i\c(y_i,u,v)$ {\bf then}\\ 
\> \> $k\gets\mbox{ a new index}$   \\
\> \> $I\gets I\cup\{k\}$ \\
\> \> $\lambda_k\gets\lambda_i-\alpha/\c(y_i,u,v)$ \\ 
\> \> $\lambda_i\gets\alpha/\c(y_i,u,v)$ \\  
\> \> $y_k\gets y_i$ \\ 
\> \> $L_k\gets L_i$ \\ 
$y_i\gets y_i+\c(y_i,u,v)(\chi_{u}-\chi_{v})$ \\
Update $L_i$ by interchanging $u$ and $v$.\\ 
$x\gets\sum_{i\in I}\lambda_iy_i$ 
\end{tabbing}}}
\end{center}
\caption{Algorithmic description of the operation $\Push(i,u,v)$.}
\label{fig:push}
\end{figure}

Each time the algorithm applies the push operation, it updates 
the set $W$ of vertices reachable from $S$ in $G(\phi)$ and the set $Z$ 
of active pairs. 
If $\Push(i,u,v)$ is nonsaturating, it makes $v$ reachable from 
$S$ in $G(\phi)$, and hence $W$ is enlarged.  Once $v$ becomes 
reachable from $S$ in $G(\phi)$, it will never become active 
again for any $i\in I$ until the algorithm finds a $\delta$-augmenting 
path or all the active pairs disappear. 
Note that we encounter at most $n$ nonsaturating pushes 
before we find a $\delta$-augmenting 
path or all the active pairs disappear. 
Each time the algorithm picks an active pair $(i,v)$ and 
applies $\Push(i,u,v)$, the vertex $v$ shifts towards the 
end of $L_i$.
Hence the algorithm picks an active pair $(i,v)$ at most $n$ times
before $v$ enters $W$ or $W\back L_i(v) = \emptyset$.  At this
point $v$ becomes inactive, and remains inactive until the next
augmentation.
Hence, for each $i\in I$ the total 
time required for processing active vertices in $L_i$ is $\O(n^2)$.

\Xomit{
If no vertices are 
added to $W$ during a $\Scan(i,u)$, it terminates with 
$L_{y_i}(u)\subseteq W$. 
}
\Lomit{
A scan continues until $W$ increases, or all arcs $(u,v)\in A(y_i)$
with $v\in V\back W$ disappear.  In the first case, the scan is
interrupted.  Thus, if a scan is completed, all pushes are saturating.
}
\Xomit{
When we perform $\Push(i,u,v)$, if the $\Push$ is nonsaturating, then 
$W$ is enlarged and a new index $k$ is introduced. 
Whether the $\Push$ is saturating or nonsaturating, 
we swap $u$ and $v$ in the given linear orderinging before the push, which then is 
a linear extension of ${\cal P}(y_i)$ for the new $y_i$. We denote 
the new linear orderinging by $L_{y_i}$ again.
If $u$ becomes not active in $L_{y_i}$, then $L_{y_i}(u)\subseteq W$ and 
$u$ will not become active again before the next flow augmentation or 
the $\delta$-scaling phase is finished. If there still exists an active
vertex $u'$ in $L_{y_i}$, we continue to scan $u'$. We can easily see that
the total number of active vertices found for $y_i$ is less than $n$.
Also, for an active vertex $u$ the number of swappings before it becomes
not active is less than $n$.
} 


We note that we could relax the 
definition of an active vertex to include any vertex $v\in V\back W$ 
whose immediate successor in $L_i$ belongs to $W$. The correctness 
argument would apply without modifications. However, care is needed
to obtain an efficient implementation.
\Lomit{
 this generic 
algorithm had to keep larger $Z$, which would require to discuss an 
implementation issue.  
 }

\Lomit{
A scan uses at most $n$ function calls, and $O(n)$ additional operations.
}

If we find a $\delta$-augmenting path, the algorithm augments 
$\delta$ units of flow along the path, which effectively increases 
$z^-(V)$ by $\delta$. We also compute an expression for $x$ as a convex 
combination of at most $n$ affinely independent extreme bases $y_i$, 
chosen from the current $y_i$'s.
This computation is 
a standard linear programming technique of transforming feasible
solutions into basic feasible solutions.  
If the set of extreme points are not 
affinely independent, there is a set of coefficients $\mu_i$ for $i\in I$ 
that is not identically zero and satisfies $\sum\mu_iy_i=0$ and 
$\sum\mu_i=0$.  
Using Gaussian elimination, we can start computing such $\mu_i$
until a dependency is detected.  At this point, we eliminate the
dependency by computing
$\theta:=\min\{\lambda_i/\mu_i\ |\ \mu_i>0\}$ and 
update $\lambda_i:=\lambda_i-\theta\mu_i$ for $i\in I$. At least one 
$i\in I$ satisfies $\lambda_i=0$. Delete such $i$ from $I$. 
We continue this procedure until we eventually obtain affine independence. 
Since a new index $k$ is added to $I$ 
only as a result of a nonsaturating push, $|I|\leq 2n$ after finding an
augmenting path.  
The bottleneck in this procedure is the
time spent computing the coefficients $\mu_i$, which is $O(n^3)$ overall.
\Lomit{
which implies that a affinely independent expression
for $x$ can be computed in
$\O(n^3)$ time.
}

A $\delta$-scaling phase ends when either $S=\emptyset$, $T=\emptyset$,
or $Z = \emptyset$.
In the last case, we have a set of vertices $W\subset V$ that
are reachable from $S$ in $G(\phi)$ such that $W\cap T = \emptyset$.
\begin{lemma}
\label{lem:tight}
If $Z=\emptyset$, then $W$ is tight for $x$. 
\end{lemma}
\begin{proof}
If $Z=\emptyset$, for each $i\in I$ the first $|W|$ vertices in $L_i$ 
must belong to $W$. 
Then it follows from (\ref{eq:greedy0}) that $y_i(W)=f(W)$. Since 
$x=\sum_{i\in I}\lambda_iy_i$ and $\sum_{i\in I}\lambda_i=1$, 
this implies $x(W)=\sum_{i\in I}\lambda_iy_i(W)=f(W)$. 
\end{proof}
\Lomit{
being disjoint with $T$ and having no leaving arcs in 
$\bigcup_{i\in I} A(y_i)$.
}

\subsection{Correctness and Complexity}

We now investigate the number of iterations in each
$\delta$-scaling phase.  To do this, we prove 
relaxed weak and strong dualities. The next lemma shows a relaxed weak
duality.
\begin{lemma} 
\label{lem:kappa}
For any base $x\in\B(f)$ and any $\delta$-feasible flow $\phi$, the 
vector $z=x-\partial\phi$ satisfies 
$z^-(V)\leq f(X)+{n\choose 2}\delta$ for any $X\subseteq V$. 
\end{lemma} 
\begin{proof}
For any $X\subseteq V$ we have $x(X)\leq f(X)$ and 
$\partial\phi(X)\geq-{n\choose 2}\delta$, and hence 
$z^-(V)\leq z(X)\leq f(X)+{n\choose 2}\delta$. 
\end{proof}

A relaxed strong duality is given as follows.
\begin{lemma}
\label{lem:delta}
At the end of each $\delta$-scaling phase, the following 
\rnp{1}--\rnp{3} hold for $x$ and $z=x-\partial\phi$. 
\begin{description}
\item[\rnp{1}] If $S=\emptyset$, then $x^-(V)\geq f(\emptyset)-n^2\delta$ and 
$z^-(V)\geq f(\emptyset)-n\delta$. 
\item[\rnp{2}] If $T=\emptyset$, then $x^-(V)\geq f(V)-n^2\delta$ 
and $z^-(V)\geq f(V)-n\delta$. 
\item[\rnp{3}] If $W$ is tight for $x$, then 
$x^-(V)\geq f(W)-n^2\delta$ and 
$z^-(V)\geq f(W)-n\delta$. 
\end{description} 
\end{lemma}
\begin{proof} 
When the $\delta$-scaling phase finishes with $S=\emptyset$, 
we have $x(v)>\partial\phi(v)-\delta\ge -n\delta$ for every 
$v\in V$, which implies $x^-(V)\ge f(\emptyset)-n^2\delta$ as well as 
$z^-(V)\geq f(\emptyset)-n\delta$. 
Similarly, when the $\delta$-scaling phase finishes with 
$T=\emptyset$, we have $x(v)<\partial\phi(v)+\delta\le n\delta$ 
for every $v\in V$, which implies 
$x^-(V)\ge x(V)-n^2\delta = f(V)-n^2\delta$ as well as 
$z^-(V)\geq x(V)-n\delta$. 

When the $\delta$-scaling phase ends with $x(W) = f(W)$ 
due to Lemma 3.1, then  
$S\subseteq W \subseteq V\back T$ and $\partial\phi(W)<0$.
\Lomit{
$S\neq\emptyset$ and 
$T\neq\emptyset$, we have a vertex subset $W\subseteq V$ such that 
$S\subseteq W\subseteq V\back T$ and there exists no arc leaving $W$ 
in $G(\phi)$ nor in $H(y_i)$ for any $i\in I$. Then we have 
$\partial\phi(W)<0$ and $y_i(W)=f(W)$ for every $i\in I$. 
Since $x\in\B(f)$ is the convex combination of $y_i$'s, the latter
implies $x(W)=f(W)$. 
}
By the definitions of $S$ and $T$, we also have
$x(v)>\partial\phi(v)-\delta\geq-n\delta$ for every $v\in V\back W$ and
$x(v)<\partial\phi(v)+\delta\leq n\delta$ for every $v\in W$.
Therefore we have $x^-(V)=x^-(W)+x^-(V\back W)\ge 
x(W)-n\delta |W|-n\delta |V\back W| = f(W)-n^2\delta$ 
as well as $z^-(V)=z^-(W)+z^-(V\back W)\ge 
x(W)-\partial\phi(W)-|W|\delta-\delta|V\back W|\geq f(W)-n\delta$. 
\end{proof}

Lemma~\ref{lem:delta} implies that at the beginning of the 
$\delta$-scaling phase, after $\delta$ is cut in half, 
$z^-(V)$ is at least $f(X)-2n\delta$ 
for some $X\subseteq V$. Making the current flow $\delta$-feasible 
decreases $z^-(V)$ by at most ${n\choose 2}\delta$. Each 
$\delta$-augmentation increases $z^-(V)$ by $\delta$. Since 
$z^-(V)$ is at most $f(X)+{n\choose 2}\delta$ 
at the end of a $\delta$-phase by Lemma~\ref{lem:kappa} 
the number of $\delta$-augmentations per phase is at most $n^2+n$ for
all phases after the first. 
Since $z^-(V)=x^-(V) \geq-nM$ at the start of the algorithm,
setting the initial $\delta = M$ is more than sufficient obtain
a similar bound on the number of augmentations in the first phase.

As an immediate consequence of Lemmas~\ref{lem:Edm} and \ref{lem:delta}, 
we also obtain the following. 
\begin{theorem}
\label{th:opt}
The algorithm obtains a minimizer of $f$ at the end of the 
$\delta$-scaling phase with $\delta<1/n^2$. 
\end{theorem}
\begin{proof}
By Lemma~\ref{lem:delta}, the output $X$ of the algorithm satisfies 
$x^-(V)\geq f(X)-n^2\delta>f(X)-1$. 
For any $Y\subseteq V$, the weak duality in Lemma~\ref{lem:Edm} 
asserts $x^-(V)\leq f(Y)$. Thus we have $f(X)-1<f(Y)$, which implies 
by the integrality of $f$ that $X$ minimizes $f$. 
\end{proof}

\begin{theorem}
Algorithm $\SFM$ runs in $\O(n^5 \log (nM))$ time.
\end{theorem}
\begin{proof}
The algorithm starts with $\delta = M$ and ends with $\delta < 1/n^2$, so
the algorithm consists of $\O(\log (nM))$ scaling phases. 
Each scaling phase finds $\O(n^2)$ $\delta$-augmenting paths. 
To find an augmenting path, we perform at most $O(n^2)$ pushes
per extreme base.
A saturating push requires $\O(1)$ time while a nonsaturating one $\O(n)$
time. Here, note that there are less than $n$ nonsaturating pushes per
augmenting path.
Hence, the time spent in pushes per augmenting path is $\O(n^3)$.  After
each augmentation, we also update the expression 
$x=\sum_{i\in I} \lambda_i y_i$, which also takes $\O(n^3)$ time
per augmentation.
Thus the overall complexity of $\SFM$ is $\O(n^5 \log(nM))$.
\Lomit{
scan each admissible pair $(i,u)$ for
$i\in I$ and $u\in V$ at most
once, and there are $\O(n^2)$ admissible pairs, since $|I| \leq 2n$.  
In addition, we make at most $n-1$ interrupted scans.  
Thus, the total number of scans per augmentation is $\O(n^2)$.
Each scan uses at most $n-1$ push operations.
\Lomit{
When a push operation is performed, we compute the
associated exchange capacity 
(see Lemma \ref{lemma:ideal}) and construct the Hasse 
diagram for the new extreme base, which takes $\O(n^2)$ time.
}
Thus the total time for a scan is $\O(n)$.
After each scan or each non-saturating push, we must determine
a new admissible pair.  This takes $O(n)$ time by traversing
the linear ordering for $y_i$ until the first vertex in $U_i$
is encountered.
}
\end{proof}

\Xomit{
When we find a $\delta$-augmenting path,
we update the expression $x=\sum_{i\in I}\lambda_iy_i$, which is carried out in
$\O(n^3)$ time. To find a $\delta$-augmenting path requires $\O(n^2)$ scans 
of admissible pairs $(i,u)$. 
It follows from (\ref{eq:W1}) and (\ref{eq:W2}) that the total time required 
for finding admissible pairs to get a $\delta$-augmenting path is $\O(n^3)$ 
by a topological sort. Each scan of an admissible pair requires $\O(n)$ push 
operations. 
When a push operation is performed, we compute the
associated exchange capacity 
(see Lemma \ref{lemma:ideal}) and constructs the Hasse 
diagram for the new extreme base, which takes $\O(n^2)$ time.
Since the total number of pushes required for getting a $\delta$-augmenting
path is $\O(n^3)$, the overall time complexity of our algorithm is 
$\O(n^7\log nM)$. 
}

In this section, we have shown a weakly polynomial-time algorithm for 
minimizing integer-valued submodular functions. The integrality of 
a submodular function $f$ guarantees that if we have a base 
$x\in {\rm B}(f)$ and a subset $X$ of $V$ such that the duality gap 
$f(X)-x^-(V)$ is less than one, $X$ is a minimizer of $f$. Except 
for this we have not used the integrality of $f$. It follows that 
for any real-valued submodular function $f: 2^V\to\R$, if we are 
given a positive lower bound $\epsilon$ for the difference between 
the second minimum and the minimum value of $f$, the present algorithm 
works for the submodular function $(1/\epsilon)f$ and runs in 
$\O(n^5\log (nM/\epsilon))$ time, where $M$ is an upper bound on 
$|f(X)|$ among $X\subseteq V$.

\section{A Strongly Polynomial-Time Algorithm} 
\label{sec:spa}
This section presents a strongly polynomial-time algorithm for 
minimizing submodular functions using the scaling algorithm 
in Section~\ref{sec:sa}. The new algorithm exploits the 
following proximity lemma. 

\begin{lemma}
\label{lem:prox}
At the end of the $\delta$-scaling phase, if $x(w)<-n^2\delta$, 
then $w$ belongs to every minimizer of $f$. 
\end{lemma} 
\begin{proof}
Let $X$ be any minimizer of $f$. 
There exists a vector $y\in\B(f)$ with 
$x^-\leq y^-$ such that $y^-(V)=f(X)$. Note that $y(v)\geq 0$ 
for each $v\in V\back X$. By Lemma~\ref{lem:delta}, there exists a 
subset $Y\subseteq V$ such that $x^-(V)\geq f(Y)-n^2\delta$. 
Then we have $y^-(w)-x^-(w)\leq y^-(V)-x^-(V)\leq f(X)-f(Y)+n^2\delta
\leq n^2\delta$. This implies $y(w)<0$ due to the assumption, 
and hence $w\in X$.  
\end{proof}

Let $f:2^V\to\R$ be a submodular function and $x\in\B(f)$ an extreme 
base whose components are bounded from above by $\eta>0$. Assume that 
there exists a subset $Y\subseteq V$ such that $f(Y)\leq -\kappa$ for 
some positive parameter $\kappa$, which will be specified later 
as $\eta/2$. We then apply the scaling algorithm 
starting with $\delta=\eta$ and the extreme base $x\in \B(f)$. After 
$\lceil\log_2(n^3\eta/\kappa)\rceil$ scaling phases, $\delta$ becomes 
less than $\kappa/n^3$. Since $x(Y)\leq f(Y)\leq -\kappa$, at least one 
element $w\in Y$ satisfies $x(w)<-n^2\delta$. 
By Lemma~\ref{lem:prox}, such an element $w$ belongs to every 
minimizer of $f$. We denote this procedure by $\Fix(f,x,\eta)$. 

We now discuss how to apply this procedure to design a strongly 
polynomial-time algorithm for minimizing a submodular function $f$. 
If $f(V)>0$, we replace the value $f(V)$ by zero. 
The set of minimizers remains the same unless the minimum value 
is zero, in which case we may assert that $\emptyset$ minimizes $f$. 

An ordered pair $(u,v)$ of distinct vertices $u,v\in V$ is said to be 
{\em compatible} with $f$ if $u\in X$ implies $v\in X$ for 
every minimizer $X$ of $f$. Our algorithm keeps 
a directed acyclic graph $D=(V,F)$ whose arcs are compatible with 
$f$. Initially, the arc set $F$ is empty. Each time the algorithm 
finds a compatible pair $(u,v)$ with $f$, it adds $(u,v)$ to $F$. 
When this gives rise to a cycle in $D$, the algorithm contracts the 
strongly connected component $U\subseteq V$ to a single vertex and 
modifies the submodular function $f$ by regarding $U$ as a singleton. 

For each $v\in V$, let $R(v)$ denote the set of vertices 
reachable from $v$ in $D$ and $f_v$ the submodular function on 
the subsets of $V\back R(v)$ defined by 
$$f_v(X)=f(X\cup R(v))-f(R(v))\quad\quad (X\subseteq V\back R(v)).$$ 
A linear ordering $(v_1,\cdots,v_n)$ of $V$ is called {\it consistent} 
with $D$ if $i<j$ implies $(v_i,v_j)\notin F$. Consider an extreme 
base $x\in\B(f)$ generated by a linear ordering $(v_1,v_2,\cdots,v_n)$
           consistent with $D$.
The extreme base generated by a consistent linear ordering is also
called {\it consistent}. 
For any consistent extreme base $x\in\B(f)$, the greedy algorithm
defines $x(v)$ by
$x(v) = f(U)-f(U\back\{v\})$ for some
$U\subseteq V$ with $R(v)\subseteq U$.  It then follows from
the submodularity of $f$ that $x$ satisfies 
$x(v)\leq f(R(v))-f(R(v)\back\{v\})$ for each $v\in V$.

In each iteration, the algorithm computes 
\begin{equation}
\label{eq:eta}
\eta=\max\{f(R(v))-f(R(v)\back\{v\})\mid v\in V\}. 
\end{equation}
If $\eta\leq 0$, then an extreme base $x\in\B(f)$ consistent with $D$ 
satisfies $x(v)\leq 0$ for each $v\in V$.  In this case
$x^-(V) = x(V) = f(V)$, which implies that 
$V$ minimizes $f$ by Lemma~\ref{lem:Edm}. 
If in addition $f(V)=0$, then the original function 
may have had a positive value of $f(V)$. Therefore, the algorithm 
returns $\emptyset$ or $V$ as a minimizer, according to whether 
$f(V)=0$ or $f(V)<0$. 

If $\eta>0$, let $u$ be an element that 
attains the maximum in the right-hand side of (\ref{eq:eta}). Then we 
have $f(R(u))=f(R(u)\back\{u\})+\eta$, which implies either 
$f(R(u))\geq\eta/2>0$ or $f(R(u)\back\{u\})<-\eta/2<0$ holds. 

In the former case ($f(R(u))\geq\eta/2$), we have 
$f_u(V\back R(u))=f(V)-f(R(u))\leq-\eta/2$. 
The algorithm finds a consistent extreme base $x\in\B(f_u)$ generated 
by a linear ordering $(v_1,\cdots,v_k)$ consistent with $D$, where 
$k=|V\back R(u)|$.  
That is,
let $x(v_1) = f_u(\{v_1\})$ and
$x(v_j)=f_u(\{v_1,v_2,\ldots, v_j\})-f_u(\{v_1,v_2,\ldots,v_{j-1}\})$
for $j=2,\ldots,k$.
Then the extreme base $x\in\B(f_u)$ satisfies 
$x(v)\leq f(R(v))-f(R(v)\back\{v\})\leq\eta$. Thus we may apply the 
procedure $\Fix(f_u,x,\eta)$ to find an element $w\in V\back R(u)$ 
that belongs to every minimizer of $f_u$. Since 
$\kappa=\eta/2$, the procedure terminates within $\O(\log n)$ 
scaling phases.  Consequently, we obtain a new pair $(u,w)$ that is 
compatible with $f$. Hence the algorithm adds the arc $(u,w)$ to $F$. 

In the latter case ($f(R(u)\back\{u\})<-\eta/2$), we compute 
an extreme base $x\in\B(f)$ consistent with $D$ by the greedy 
algorithm, and then apply the procedure $\Fix(f,x,\eta)$ to find 
an element $w\in R(u)$ that belongs to every minimizer of $f$. 
Since $x(v)\leq\eta$ for every $v\in V$ and $\kappa=\eta/2$,
the procedure terminates within $\O(\log n)$ scaling phases. 
Note that every minimizer of $f$ includes $R(w)$. Thus it suffices 
to minimize the submodular function $f_w$, which is now defined on 
a smaller underlying set.    
Figure \ref{fig:spta} provides a formal 
description of the strongly polynomial-time algorithm. 

\begin{figure}
\framebox[\textwidth]{\parbox{6in}{
\begin{tabbing}
****\=***\=***\=***\=***\= \hspace{2.5in} \=  \kill
{\tt Input:}\ $f:2^V\rightarrow \R$ \\
{\tt Output:} $X\subseteq V$ minimizing $f$ \\ 
\> \\
{\bf Initialization:} \\
\> $X\gets\emptyset$ \\  
\> $F\gets\emptyset$ \\
{\bf While} $V\neq\emptyset$ {\bf do} \\
\> {\bf If} $f(V)>0$ {\bf then} $f(V)\gets 0$ \\ 
\> $\eta\gets\max\{f(R(v))-f(R(v)\back\{v\})\mid v\in V\}$ \\
\> {\bf If} $\eta\leq 0$ {\bf then} {\bf break} \\ 
\> Let $u\in V$ attain the maximum above. \\
\> {\bf If} $f(R(u))\geq\eta/2$ {\bf then} \\ 
\> \> Find a consistent extreme base $x\in\B(f_u)$ by the greedy algorithm. \\
\> \> $w\gets\Fix(f_u,x,\eta)$ \\ 
\> \> {\bf If} $u\in R(w)$ {\bf then} \\
\> \> \> Contract $\{v\mid v\in R(w),\,u\in R(v)\}$ to a single vertex. \\
\> \> {\bf Else} $F\gets F\cup\{(u,w)\}$ \\ 
\> {\bf Else} \\
\> \> Find a consistent extreme base $x\in\B(f)$ by the greedy algorithm. \\
\> \> $w\gets\Fix(f,x,\eta)$ \\ 
\> \> $V\gets V\back R(w)$ \\ 
\> \> $f\gets f_w$ \\ 
\> \> Find a subset $Q$ of the original underlying set 
represented by $R(w)$. \\
\> \> $X\gets X\cup Q$ \\ 
{\bf If} $f(V)<0$ {\bf then} \\ 
\> \> Find a subset $Q$ of the original underlying set 
represented by $V$. \\
\> \> $X\gets X\cup Q$ \\ 
{\bf End}.
\end{tabbing}
}}
\caption{A strongly polynomial-time algorithm for 
submodular function minimization.}
\label{fig:spta}
\end{figure}

\begin{theorem}
The algorithm in {\em Figure~\ref{fig:spta}} computes the minimizer of
a submodular function in $\O(n^7\log n)$ time, which is strongly polynomial.
\end{theorem}
\begin{proof}
Each time we call the procedure $\Fix$, the algorithm adds a new arc to $D$ or 
deletes a set of vertices. This can happen at most $n^2$ times. Thus the 
overall running time of the algorithm is $\O(n^7\log n)$, which is 
strongly polynomial. 
\end{proof}

\section{Concluding Remarks}
\label{sec:cr} 
This paper presents a strongly polynomial-time algorithm for 
minimizing submodular functions defined on Boolean lattices. We now 
briefly discuss minimizing submodular functions defined on more general 
lattices. 

Consider a submodular function $f: {\cal D}\to\R$ defined 
on a distributive lattice ${\cal D}$ represented by a poset $\P$ on
$V$. Then the associated base polyhedron is unbounded in general (see \cite{Fuji91}).

An easy way to minimize such a function $f$ is to consider the 
reduction of $f$ by a sufficiently large vector. As described in 
\cite[p.~56]{Fuji91}, we can compute an upper bound $\hat{M}$ 
on $|f(X)|$ $(X\in {\cal D})$. Let $f'$ be the rank function 
of the reduction by a vector with each component being equal to 
$\hat{M}$. The submodular function $f'$ is defined on $2^V$ and 
the set of minimizers of $f'$ coincides with that of $f$. Thus, 
we may apply our algorithms. However, each evaluation of the 
function value of $f'$ requires $\O(n^2)$ elementary operations 
in addition to a single call for the evaluation of $f$. Consequently, 
this approach takes $\O(n^7\min\{\log(n\hat{M}),n^2\log n\})$ time.  

Alternatively, we can slightly extend the algorithms in 
Sections \ref{sec:sa} and \ref{sec:spa} by keeping the base 
$x\in\B(f)$ as a convex combination of extreme bases $y_i$'s 
plus a vector in the characteristic cone of $\B(f)$. The latter 
can be represented as a boundary of a nonnegative flow in the Hasse 
diagram of $\P$. This extension enables us to minimize $f$ in 
$\O(n^5\min\{\log(n\hat{M}),n^2\log n\})$ time. 

Submodular functions defined on modular lattices naturally arise 
in linear algebra. Minimization of such functions has a significant 
application to computing canonical forms of partitioned matrices 
\cite{Ito+Iwata+Murota94,Iwata+Murota95}. It remains an interesting
open problem to develop an efficient algorithm for minimizing 
submodular functions on modular lattices, even for those specific 
functions that arise from partitioned matrices. 

Independently of this work, and almost simultaneously, Schrijver has 
also developed a combinatorial, strongly polynomial-time
algorithm for submodular function minimization~\cite{Schrijver99}.  
His algorithm
also extends Cunningham's approach.  However, the resulting algorithm
is quite different from ours.

\section*{Acknowledgments}

We are grateful to Bill Cunningham, Michel Goemans, and Maiko Shigeno
for their useful comments. 


\begin{thebibliography}{99}

\bibitem{Bixby+Cunningham+Topkis} 
{R. E. Bixby}, {W. H. Cunningham}, and {D. M. Topkis}:
Partial order of a polymatroid extreme point, 
{\it Math. Oper. Res.}, {\bf 10} (1985), 367--378.

 
\bibitem{Cunningham84} {W. H. Cunningham}: 
Testing membership in matroid polyhedra, 
{\it J. Combinatorial Theory}, {\bf B36} (1984), 161--188.

\bibitem{Cunningham85} {W. H. Cunningham}: 
On submodular function minimization, 
{\it Combinatorica}, {\bf 5} (1985), 185--192. 

\bibitem{Edm70} 
{J. Edmonds}: Submodular functions, matroids, and certain polyhedra,
{\it Combinatorial Structures and Their Applications}, 
R. Guy, H. Hanani, N. Sauer, and J. Sch\"onheim, eds., 
Gordon and Breach, 69--87, 1970. 

\bibitem{Edmonds+Giles77} {J. Edmonds} and {R. Giles}: 
A min-max relation for submodular function on graphs, 
{\it Ann. Discrete Math.}, {\bf 1} (1977), 185--204.

\bibitem{Edmonds+Karp72} {J. Edmonds} and {R. Karp}: 
Theoretical improvements in algorithmic efficiency for network 
flow problems, {\it J. ACM}, {\bf 19} (1972), 248--264. 

\bibitem{FIM99} 
{L. Fleischer}, {S. Iwata}, and {S. T. McCormick}: 
A faster capacity scaling algorithm for submodular flow, 1999. 

\bibitem{Frank+Tardos87} {A. Frank} and {\'E. Tardos}: 
An application of simultaneous Diophantine approximation in 
combinatorial optimization, {\it Combinatorica}, {\bf 7} (1987),
49--65. 

\bibitem{Fuji78} {S. Fujishige}: 
Polymatroidal dependence structure of a set of random variables, 
{\it Information and Control}, {\bf 39} (1978), 55--72.  

\bibitem{Fuji80} {S. Fujishige}: 
Lexicographically optimal base of a polymatroid with respect to a 
weight vector, {\it Math. Oper. Res.}, {\bf 5} (1980), 186--196. 

\bibitem{Fuji84} {S. Fujishige}: 
Submodular systems and related topics, 
{\it Math. Programming Study}, {\bf 22} (1984), 113--131. 

\bibitem{Fuji91} {S. Fujishige}: 
{\it Submodular Functions and Optimization}, 
North-Holland, 1991. 

\bibitem{Goemans+Ramakrishnan95}
{M. X. Goemans} and {V. S. Ramakrishnan}: 
Minimizing submodular functions over families of subsets, 
{\it Combinatorica}, {\bf 15} (1995), 499--513. 

\bibitem{GLS81} 
{M. Gr\"{o}tschel}, {L. Lov\'{a}sz}, and {A. Schrijver}:
The ellipsoid method and its consequences in combinatorial optimization,
{\it Combinatorica}, {\bf 1} (1981), 169--197.

\bibitem{GLS88} 
{M. Gr\"{o}tschel}, {L. Lov\'{a}sz}, and {A. Schrijver}:
{\it Geometric Algorithms and Combinatorial Optimization}, 
Springer-Verlag, 1988.

\bibitem{Han79} {T.-S. Han}: 
The capacity region of general multiple-access channel with 
correlated sources, {\it Information and Control}, {\bf 40} (1979),
37--60.
 
\bibitem{Hoppe+Tardos95} {B. Hoppe} and {\'E. Tardos}: 
The quickest transshipment problem, 
{\it Proceedings of 5th ACM/SIAM Symposium on Discrete Algorithms}
(1995), 512--521. 

\bibitem{Ito+Iwata+Murota94} 
{H. Ito}, {S. Iwata}, and {K. Murota}: 
Block-triangularization of partitioned matrices under 
similarity/equivalence transformations, 
{\it SIAM J. Matrix Anal. Appl.}, {\bf 15} (1994), 1226--1255.  

\bibitem{Iwata97} {S. Iwata}: 
A capacity scaling algorithm for convex cost submodular flows, 
{\it Math. Programming}, {\bf 76} (1997), 299--308.

\bibitem{Iwata+Murota95} {S. Iwata} and {K. Murota}: 
A minimax theorem and a Dulmage-Mendelsohn type decomposition for 
a class of generic partitioned matrices, 
{\it SIAM J. Matrix Anal. Appl.}, {\bf 16} (1995), 719--734.  

\bibitem{Lovasz83} {L. Lov\'asz}: 
Submodular functions and convexity.  
{\it Mathematical Programming --- The State of the Art},  
A.~Bachem, M.~Gr\"otschel and B.~Korte, eds., 
Springer-Verlag, 1983, 235--257. 


\bibitem{Nagamochi+Ibaraki92} {H. Nagamochi} and {T. Ibaraki}: 
Computing edge-connectivity in multigraphs and capacitated graphs,  
{\it SIAM J. Discrete Math.}, {\bf 5} (1992), 54--64.


\bibitem{Narayanan95} {H. Narayanan}: 
A rounding technique for the polymatroid membership problem, 
{\it Linear Algebra Appl.}, {\bf 221} (1995), 41--57.

\bibitem{Queyranne98} {M. Queyranne}: 
Minimizing symmetric submodular functions, 
{\it Math. Programming}, {\bf 82} (1998), 3--12.


\bibitem{Schrijver99} A. Schrijver: A combinatorial algorithm 
minimizing submodular functions in strongly polynomial time, 1999. 


\bibitem{Sohoni92} {M. A. Sohoni}: 
Membership in submodular and other polyhedra.
Technical Report TR-102-92, Department of Computer Science and Engineering,
Indian Institute of Technology, Bombay, India, 1992. 

\bibitem{Tamir93} {A. Tamir}: 
  A unifying location model on tree graphs based on submodularity
 properties, {\it Discrete Appl. Math.}, {\bf 47} (1993),
 275--283. 



\end{thebibliography}
\end{document}